\documentclass [11pt]{amsart}
\usepackage {amsmath, amssymb, amscd, mathrsfs, hyperref, graphicx, color, slashed, enumerate, url}
\usepackage[text={6in,8.5in},centering,letterpaper,dvips]{geometry}
\usepackage[all,cmtip]{xy}

\newtheorem {theorem}{Theorem}[section]

\newtheorem {definition}[theorem]{Definition}

\theoremstyle{remark}
\newtheorem {remark}[theorem]{Remark}
\newtheorem {example}[theorem]{Example}

\DeclareFontFamily{U}{mathx}{\hyphenchar\font45}
\DeclareFontShape{U}{mathx}{m}{n}{
      <5> <6> <7> <8> <9> <10>
      <10.95> <12> <14.4> <17.28> <20.74> <24.88>
      mathx10
      }{}
\DeclareSymbolFont{mathx}{U}{mathx}{m}{n}
\DeclareFontSubstitution{U}{mathx}{m}{n}
\DeclareMathAccent{\widecheck}{0}{mathx}{"71}

\def\polhk#1{\setbox0=\hbox{#1}{\ooalign{\hidewidth
    \lower1.5ex\hbox{`}\hidewidth\crcr\unhbox0}}}  
\def\Z {{\mathbb{Z}}}
\def\R {{\mathbb{R}}}
\def\C {{\mathbb{C}}}
\def\Q {{\mathbb{Q}}}
\def\rr{\R}
\def\Conf {\mathcal{C}}
\def\I{\mathfrak{I}}
\def\Cat{\mathfrak{C}}
\def\HTop{\mathrm{HTop}}

\def\Spin {\mathbb{S}}
\def\dirac {\slashed{\partial}}
\def\Dirac{\slashed{D}}
\def\tH{\tilde{H}}
\def\rp{\mathbb{RP}}
\def\cp{\mathbb{CP}}
\def\hp{\mathbb{HP}}
\def\F {\mathbb{F}_2}
\def\del {\partial}
\def\f {\mathbb{F}}
\def\SW{\mathit{SW}}
\def\S{\mathcal{S}}
\def\t{\mathfrak{t}}

\DeclareMathOperator{\id}{\operatorname{id}}
\def\ind{\operatorname{ind}}

\def\Sq{\operatorname{Sq}}

\def\Inv{\operatorname{Inv}}
\def\Int {\operatorname{int}}
\def\To {\longrightarrow}
\def\su{\mathfrak{su}}
\def\spinc{\operatorname{Spin}^c}
\def\sutwo {\operatorname{SU}(2)}
\def\tM{\widetilde{M}}

\def\Cat{\mathfrak{S}}

\def\sl{\mathfrak{sl}}

\def\link{\mathit{link}}

\def\swf{\operatorname{SWF}}
\def\swfh{\mathit{SWFH}}
\def\swfc{\mathit{SWFC}}

\def\pin {\operatorname{Pin}(2)}

\newcommand{\HMto}{\widecheck{\mathit{HM}}}

\def\H {\mathbb{H}}

\def\csd{\mathit{CSD}}

\def\G {\mathcal{G}}
\newcommand{\s}{\mathfrak{s}}

\def\pt {\operatorname{pt}}

\begin{document}

\title[The Conley index, gauge theory, and triangulations]{The Conley index, gauge theory, and triangulations}

\author[Ciprian Manolescu]{Ciprian Manolescu}
\thanks {The author was supported by NSF grant DMS-1104406.}
\address {Department of Mathematics, UCLA, 520 Portola Plaza\\ Los Angeles, CA 90095}
\email {cm@math.ucla.edu}

\begin{abstract}
This is an expository paper about Seiberg-Witten Floer stable homotopy types. We outline their  construction, which is based on the Conley index and finite dimensional approximation. We then describe several applications, including the disproof of the high-dimensional triangulation conjecture.    
\end {abstract}

\maketitle

\section{Introduction}

The Conley index is an important topological tool in the study of dynamical systems.  Conley's monograph \cite{ConleyBook} is the standard reference on this subject; see also \cite{SalamonConley, MischaikowSurvey, MischaikowMrozek} for more recent expositions. In symplectic geometry, the Conley index was notably used to prove the Arnol'd conjecture for the $n$-dimensional torus \cite{ConleyZehnder}. Furthermore, it inspired the development of Floer homology \cite{FloerWitten, FloerLagrangian, SalamonSurvey}, which is an infinite dimensional variant of Morse theory. Apart from symplectic geometry, Floer homology appears in the context of gauge theory, where it produces three-manifold invariants starting from either the instanton (Yang-Mills) or the monopole (Seiberg-Witten) equations \cite{FloerInstanton, MarcolliWang, KMbook, FroyshovSW}.

In \cite{Spectrum}, the author used the Conley index more directly to define a version of ($S^1$-equivariant) Seiberg-Witten Floer homology. The strategy was to approximate the Seiberg-Witten equations by a gradient flow in finite dimensions, and then take the homology of the appropriate Conley index. This is very similar in spirit to the work of G{\polhk{e}}ba, Izydorek and Pruszko \cite{GIP}, who defined a Conley index for flows on Hilbert spaces using finite dimensional approximation. However, the Seiberg-Witten case is more difficult analytically, because there is no monopole map from a single Hilbert space to itself; rather, the monopole map takes a Sobolev space to one of lower regularity.

Compared to the more traditional constructions of monopole Floer homology \cite{KMbook, MarcolliWang, FroyshovSW}, the method in \cite{Spectrum} has the following advantages:
\begin{itemize}
\item It avoids dealing with transversality issues, such as finding generic perturbations: To take the Conley index, one does not have to ensure that the gradient flow is Morse-Smale;
\item It makes it easy to incorporate symmetries of the equations: the $S^1$-symmetry in \cite{Spectrum}, the $\pin$-symmetry used in \cite{beta}, as well as finite group symmetries coming from  coverings of three-manifolds \cite{LidmanM};
\item It yields more than just Floer homologies: The equivariant stable homotopy type of the Conley index is a three-manifold invariant. One can then apply to it other generalized homology functors, such as equivariant K-theory \cite{kg}.  
\end{itemize}

In contrast to the Kronheimer-Mrowka construction of monopole Floer homology \cite{KMbook} (which works for all three-manifolds), one limitation of the Conley index approach is that so far it has only been developed for manifolds with $b_1=0$ \cite{Spectrum}, and partially for manifolds with $b_1=1$ \cite{PFP}. The main difficulty is that one needs to find suitable finite dimensional approximations. This is easy to do for rational homology spheres, when the configuration space is a Hilbert space---the approximations are given by finite dimensional subspaces. However, it becomes harder for higher $b_1$, due to the presence of an algebraic-topologic obstruction called the polarization class; we refer to \cite[Appendix A]{PFP} for more details.

This article is meant as an introduction to the finite dimensional approximation / Conley index technique in Seiberg-Witten Floer theory. We discuss several consequences, and in particular highlight the following application of $\pin$-equivariant Seiberg-Witten Floer homology:

\begin{theorem}[\cite{beta}]
\label{thm:tc}
There exist non-triangulable $n$-dimensional topological manifolds for every $n \geq 5$.
\end{theorem}

Previously, non-triangulable manifolds have been shown to exist in dimension four by Casson \cite{Casson}. The proof of Theorem~\ref{thm:tc} rests on previous work of Galewski-Stern and Matumoto \cite{GS, Matumoto}, who reduced the problem to a question about homology cobordism in three dimensions. That question can then be answered using Floer theory.

The paper is organized as follows. In Section~\ref{sec:Conley} we give an overview of Conley index theory, focusing on gradient flows and the relation to Morse theory. In Section~\ref{sec:SWF} we describe the construction of the Seiberg-Witten Floer stable homotopy type for rational homology three-spheres. In Section~\ref{sec:tc} we present the historical background to the triangulation problem, and sketch its solution. Finally, in Section~\ref{sec:other} we discuss other topological applications.

\medskip

\noindent {\bf Acknowledgements.} The author is indebted to Mike Freedman, Rob Kirby, Peter Kronheimer, Frank Quinn, Danny Ruberman and Ron Stern for helpful conversations related to the triangulation problem. Comments and suggestions by Rob Kirby, Mayer Landau, Tye Lidman and Frank Quinn on a previous draft have greatly improved this article.

\section{The Conley index}
\label{sec:Conley}

\subsection{Morse complexes}
\label{sec:morse}
Let $M$ be a closed Riemannian manifold. Given a Morse-Smale function $f: M \to \rr$, there is an associated Morse complex $C_*(M, f)$. The generators are the critical points of $f$ and the differential is given by
\begin{equation}
\label{eq:del}
 \del x = \sum_y n_{xy} y,
 \end{equation}
where $n_{xy}$ is the signed count of index $1$ gradient flow lines between $x$ and $y$. The Morse homology $H_*(M, f)$ is isomorphic to the usual singular homology of $M$.

Let us investigate what happens if we drop the compactness assumption. In general, it may no longer be the case that $\del^2 = 0$:

\begin{example}
\label{ex:1}
Suppose we have a Morse-Smale function $f$ on a surface, a gradient flow line from a local maximum $x$ to a saddle point $y$, and another flow line from $y$ to a local minimum $z$. Let $M$ be a small open neighborhood of the union of these two flow lines. Then the restriction of $f$ to $M$ does not yield a Morse chain complex: we have $\del^2 x = \pm \del y = \pm z$.
\end{example}

In order to obtain a Morse complex on a non-compact manifold $M$, we need to impose an additional condition. Let $f: M \to \rr$ be Morse-Smale. Some gradient flow lines of $f$ connect critical points, while others may escape to the ends of the manifold, in positive and/or in negative time (and either in finite or in infinite time). Let us denote by $\S \subseteq M$ the subset of all points that lie on flow lines connecting critical points. (In particular, $S$ includes all the critical points.) In Example~\ref{ex:1}, the set $\S$ is not compact. This is related to the non-vanishing of $\del^2$: The moduli space of flow lines from $x$ to $z$ is a one-dimensional manifold $M \cong (0,1)$, and the broken flow line through $y$ only gives a partial compactification of $M$ (by a single point). 

Therefore, let us assume that $\S$ is compact. Then, the same proof as in the case when $M$ is compact shows that the differential $\del$ given by \eqref{eq:del} satisfies $\del^2=0$. We obtain a Morse complex $C_*(M, f)$. The next question is, what does the Morse homology $H_*(M, f)$ compute in this case? As the reader can check in simple examples, it does not give the singular homology of either $M$ or $\S$. The answer turns out to be the homology of the Conley index of $\S$, which we now proceed to define.

\subsection{The Conley index}
\label{sec:ci}
Although in this paper we will only need the Conley index in the setting of gradient flows, let us define it  more generally. Following \cite{ConleyBook}, suppose that we have a one-parameter subgroup $\varphi= \{\varphi_t \}$ of diffeomorphisms of an $n$-dimensional manifold $M$, and a compact subset $N \subseteq M$. Let 
\[
\Inv (N, \varphi) = \{x \in N \mid \varphi_t (x) \in N \text{ for all } t \in \R \}.
\]

A compact subset $N \subseteq M$ is called an {\em isolating neighborhood} if $\Inv(N, \varphi) \subseteq \Int N$. We also define an  {\em isolated invariant set} to be a subset $\S \subseteq M$ such that $\S = \Inv(N, \varphi)$ for some isolating neighborhood $N$. Note that isolated invariant sets are compact.

\begin{definition}
\label{def:indexpair}
Let $\S$ be  an isolated invariant set. An {\em index pair} $(N,L)$ for $\S$ is a pair of compact sets $L \subseteq N \subseteq M$ such that:
\begin{enumerate}[(i)]
\item $\Inv (N - L, \varphi) = \S \subset \text{int }(N - L)$.
\item $L$ is an {\em exit set} for $N$; that is, for all $x \in N$, if there exists $t >0$ such that $\varphi_t(x)$ is not in $N$, then there exists $0 \leq \tau < t$ with $\varphi_\tau(x) \in L$. 
\item $L$ is {\em positively invariant} in $N$; that is, if $x \in L$ and $t>0$ are such that $\varphi_s(x) \in N$ for all $0 \leq s \leq t$, then $\varphi_s(x)$ is in $L$ for $0 \leq s \leq t$. 
\end{enumerate}
\end{definition}

It was proved by Conley \cite{ConleyBook} that any isolated invariant set $\S$ admits an index pair. The {\em Conley index} for an isolated invariant set $\S$  is defined to be the based homotopy type
$$I(\varphi, \S) :=(N/L,[L]).$$  

\begin{theorem}[\cite{ConleyBook}]
$(a)$ The Conley index $I(\varphi, \S)$ is an invariant of the triple $(M,\varphi,\S)$. 

$(b)$ The Conley index is invariant under continuation: If we have a smooth family of flows $\varphi^\lambda = \{\varphi^{\lambda}_t\}, \lambda \in [0,1]$, and $N$ is an isolating neighborhood in every $\varphi^{\lambda}$, then the Conley index for $\S_\lambda = \Inv(N, \varphi_\lambda)$ in the flow $\varphi^{\lambda}$ is independent of $\lambda$.
\end{theorem}

\begin{example}
Suppose $\varphi$ is the downward gradient flow of a Morse function, and $\S = \{x\}$ consists of a single critical point of Morse index $k$. We can find an isolating neighborhood $N$ for $\{x\}$ of the form $D^k  \times D^{n-k}$, with $L = \del D^k \times D^{n-k}$ being the exit set. We deduce that the Conley index of $\{x\}$ is the homotopy type of $S^k$, so $k$ can be recovered from $I(\varphi, \S)$. Thus, we can view the Conley index as a generalization of the usual Morse index.
\end{example}

\begin{figure}
\begin{center}
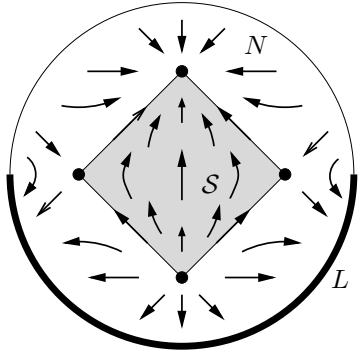
\end{center}
\caption {{\bf An isolating invariant set in a gradient flow.} The set $\S$ is shaded. The disk $N$ is an isolating neighborhood for $\S$, and $L \subset \del N$ is the exit set.} 
\label{fig:flow}
\end{figure}

In practice, it is helpful to know that we can find index pairs with certain nice properties. 
For any isolated invariant set $\S$, we can choose an index pair $(N, L)$ such that $N$ and $L$ are finite CW complexes. In fact, more is true: We can arrange so that $N$ is an $n$-dimensional manifold with boundary, and $L \subset \del N$ is an $(n-1)$-dimensional manifold with boundary. See Figure~\ref{fig:flow}. This is useful, for example, when relating the Conley index in a forward flow $\varphi$ to the Conley index in the reverse flow $\bar \varphi$, given by $\bar \varphi_t = \varphi_{-t}$. Then $\S$ is also an isolated invariant set for $\bar \varphi$. Furthermore, we can arrange so that an index pair for $\S$ in $\bar \varphi$ is $(N, L')$, where $L' \subset \del N$ is the closure of $(\del N) -L$. From here we get that, if the ambient manifold $M$ is a vector space, then $I(\varphi, \S)$ and $I(\bar \varphi, \S)$ are Spanier-Whitehead dual with respect to $M$; see \cite{McCord, Cornea} for details. 

\begin{remark}
\label{rem:css}
Although we have defined the Conley index as a homotopy type, something stronger is true. Given two different choices of index pair $(N_1, L_1)$ and $(N_2, L_2)$ for the same $\S$, they are related by an equivalence whose homotopy class is canonical. In other words, to each $\S$ one can associate a connected simple system, i.e., a subcategory $\I=\I(\varphi, \S)$ of a given category $\Cat$ (in this case the homotopy category $\HTop*$ of pointed topological spaces), such that there is exactly one morphism between any two objects in $\I$. Having a connected simple system is sometimes rephrased by saying that we have an element in the category $\Cat$ that is well-defined up to canonical isomorphism in $\Cat$. (Note that in our case, isomorphism in $\HTop*$ means based homotopy equivalence.) 
\end{remark}

Let us now return to the case of gradient flows considered in Section~\ref{sec:morse}. We have:

\begin{theorem}[Floer \cite{FloerWitten}] 
Let $\S$ be an isolated invariant set for a Morse-Smale gradient flow $\varphi$. Then, the Morse homology computed from the set of all critical points and flow lines in $\S$ is isomorphic to the reduced homology of the Conley index $I(\varphi, \S)$. 
\end{theorem}

See also \cite{RotV} for an extension of this result to more general flows.

\subsection{The equivariant Conley index}
Floer \cite{FloerConley} and Pruszko \cite{Pruszko} refined Conley index theory to the equivariant setting. Precisely, let $G$ be a compact Lie group acting smoothly on a manifold $M$, preserving a flow $\varphi$ and an isolated invariant set $\S$. Then, there exists a $G$-invariant index pair $(N, L)$ for $\S$, and the Conley index 
$$I_G(\varphi, \S):=(N/L, [L])$$ 
is well-defined up to canonical $G$-equivariant homotopy equivalence. Moreover, G{\polhk{e}}ba \cite[Proposition 5.6]{Geba} showed that $I_G(\varphi, \S)$ is the based homotopy type of a finite $G$-CW complex. 

The discussion of duality for the Conley indices in the forward and reverse flow extends to the equivariant setting.

\section{Seiberg-Witten Floer homology}
\label{sec:SWF}

The Seiberg-Witten equations \cite{SW1, SW2, Witten} play a fundamental role in low-dimensional topology. They yield the Seiberg-Witten invariants of closed four-manifolds. When one cuts a closed four-manifold $W$ along a three-manifold $Y$, the invariant of $W$ can be recovered from relative invariants of the two pieces. The latter are elements in a group associated to $Y$, called the Seiberg-Witten (or monopole) Floer homology.

In this section we will not discuss four-manifolds much, but rather focus on dimension three, and on the case of rational homology spheres. We will describe various approaches to the construction of Seiberg-Witten Floer homology in this setting. In particular, we will present the Conley index method, which also gives rise to the Seiberg-Witten Floer stable homotopy type. We will give examples and discuss a few properties of the invariants.

\subsection{The Seiberg-Witten equations in dimension three} \label{sec:sw}
Let $Y$ be a closed, oriented $3$-manifold with $b_1(Y)=0$, and let $g$ be a Riemannian metric on $Y$. A $\spinc$ structure $\s$ on $Y$ consists of a rank two Hermitian vector bundle $\Spin$, together with a Clifford multiplication $\rho: TY \to \su(\Spin)$ which maps $TY$ isometrically to the space of traceless, skew-adjoint endomorphisms of $\Spin$. The multiplication $\rho$ can be extended to real $1$-forms by duality, and then complexified to give a map $\rho: T^*Y \otimes \C \to \sl(\Spin)$. There is an associated Dirac operator $\dirac : \Gamma(\Spin) \to \Gamma(\Spin)$.

We define the configuration space:
\[
\Conf(Y,\s) = i\Omega^1(Y) \oplus \Gamma(\Spin).
\]
For a pair $(a, \phi) \in \Conf(Y, \s)$, the Seiberg-Witten equations are:
\begin{equation}
\label{eq:sw}
*da + \tau(\phi,\phi) = 0, \; \; \dirac \phi + \rho(a)\phi = 0,  
\end{equation}
where $\tau(\phi,\phi) = \rho^{-1}(\phi \otimes \phi^*)_0 \in \Omega^1(Y; i\R)$ and the subscript $0$ denotes the trace-free part. They are invariant with respect to the action of the gauge group $\G = C^\infty(Y,S^1)$, which acts on $\Conf(Y,\s)$ by $u \cdot (a,\phi) = (a - u^{-1}du,u \cdot \phi)$. 

For short, we write the equations \eqref{eq:sw} as 
$$\SW(a, \phi)=0.$$ 

Note that $\Conf(Y, \s)$ is only a Fr{\'e}chet space, not a Banach space. It is helpful to consider its $L^2_k$ Sobolev completions $\Conf_k(Y, \s)$ for large $k$. Then $\SW$ can be viewed as a map from $\Conf_k(Y, \s)$ to $\Conf_{k-1}(Y,\s)$.

The Seiberg-Witten map is the formal gradient flow of a functional on $\Conf(Y, \s)$ called the Chern-Simons-Dirac ($\csd$) functional. In this context, we expect to be able to define Floer homology 
by analogy with Morse homology in finite dimensions. Since the Seiberg-Witten solutions come in orbits of $\G$, in order to get isolated critical points we should divide by the gauge action. If we first divide by the smaller group $\G_0 \subset \G$ consisting of $u=e^{i \xi}$ with $\int_Y \xi = 0$, the quotient $\Conf(Y, \s)/\G_0$ can be identified with the Coulomb slice:
$$V = i \ker d^* \oplus \Gamma(\Spin) \subset \Conf(Y, \s).$$ 

We still have a leftover action by $S^1$ given by the constant gauge transformations, $e^{i\theta} : (a, \phi) \mapsto (a, e^{i\theta} \phi)$. If we try to divide $V$ by $S^1$ we would get a singularity at the origin. Instead, it is helpful to distinguish between two types of solutions to the Seiberg-Witten equations on $V$:
\begin{enumerate}[(i)]
\item {\em reducibles}, i.e., fixed by $S^1$. There is a unique reducible solution in every $\spinc$ structure, namely $(a, \phi)=(0,0)$; 
\item {\em irreducibles}, i.e., having a free orbit under the $S^1$ action.
\end{enumerate}

Generically, we expect that there are finitely many irreducibles (modulo $S^1$). Furthermore, Seiberg-Witten Floer homology should have an $S^1$-equivariant flavor, in the form of a module over the equivariant cohomology of a point
$$H^*_{S^1}(\pt) \cong H^*(\cp^{\infty}) = \Z[U],$$
where $U$ is in degree $2$ (and hence acts on homology by lowering degree by $2$).
The Floer homology is constructed from a complex $\swfc^{S^1}_*(Y, \s, g)$ composed of one copy of $H^{S^1}_*(\pt)$ for the reducible:
\begin{equation}
\label{eq:Utail}
 \xymatrixcolsep{.7pc}
\xymatrix{
 \Z  &  0 & \Z \ar@/_1pc/[ll]_{U} & 0 & \Z \ar@/_1pc/[ll]_{U} & 0 & \dots \ar@/_1pc/[ll]_{U} 
} 
\end{equation}
and one copy of $H^{S^1}_*(S^1) \cong H_*(\pt) \cong \Z$ for each irreducible. The generators are connected by the differential (which counts gradient flow lines, and decreases grading by $1$) as well as by the action of $U$ by slant product (which counts flow lines of a certain type, and decreases degree by $2$). In particular, the ``infinite $U$-tail'' in \eqref{eq:Utail} (coming from the reducible) can interact with the irreducibles through $\del$ or $U$. 

The complex $\swfc^{S^1}_*(Y, \s, g)$ depends on the metric $g$, but its homology $\swfh_*^{S^1}(Y, \s)$ does not. Note that $\swfh_*^{S^1}(Y, \s)$ still decomposes into an infinite $U$-tail and a finite Abelian group. The $U$-tail in homology could be shorter (that is, start in a higher degree) than the one in the Floer complex if $\del$ maps some of the irreducibles to some of the reducible generators, thus cancelling them in homology. The $U$-tail in homology could also be longer (that is, start in a lower degree) if $U$ maps some of the reducible generators to irreducibles. This observation will be important when we discuss Fr{\o}yshov-type invariants in Section~\ref{sec:beta}.

In order to make the above construction of $\swfh_*^{S^1}(Y, \s)$ rigorous, one needs to find a good class of perturbations for the Seiberg-Witten equations, so that the resulting equivariant flow is Morse-Bott-Smale. This was the approach taken by Marcolli and Wang in \cite{MarcolliWang}. Nevertheless, dealing with Morse-Bott-Smale transversality directly is rather technical, and there are various ways to get around it:
\begin{itemize}
\item In \cite{KMbook}, Kronheimer and Mrowka replaced $\Conf(Y, \s)$ by its blow-up $\Conf^\sigma(Y, \s)$ consisting of triples $(a, s, \phi) \in i\Omega^1(Y) \oplus \R \oplus \Gamma(\Spin)$ with $s \geq 0$ and $\| \phi \|_{L^2} =1$. The $L^2_k$ completion of the quotient of the blow-up by $\G$ is a Hilbert manifold with boundary, and one can do Floer theory on it instead of equivariant Floer theory on $\Conf(Y, \s)$. The Kronheimer-Mrowka version of $\swfh_*^{S^1}(Y, \s)$ is denoted $\HMto(Y, \s)$; 
\item In \cite{FroyshovSW}, Fr{\o}yshov used non-exact perturbations to get rid of the reducible solution, and then took a limit of the resulting irreducible Floer groups;
\item In \cite{Spectrum}, the author used finite dimensional approximation and then took the $S^1$-equivariant homology of the corresponding Conley index. This approach will be discussed in more detail below. Roughly, working in finite dimensions allows us to avoid Morse theory entirely, and employ singular homology instead.
\end{itemize}

In all of these constructions, the result is a Seiberg-Witten Floer homology that can be shown to be of the form advertised above (a $\Z[U]$-module, with an infinite $U$-tail and a part that is finitely generated over $\Z$). The infinite $U$-tail can be defined intrinsically in terms of the module $\swfh_*^{S^1}(Y, \s)$, as the intersection of the images of $U^k$ over all $k \geq 0$.

\subsection{Finite-dimensional approximation}
We now sketch the construction in \cite{Spectrum}. When reduced to the Coulomb slice $V$, the Seiberg-Witten map $\SW$ can be written as a sum $$\ell + c: V \to V,$$ where $\ell=(*d, \dirac)$ is the linearization of $\SW$ at $(0,0)$. The map $\ell$ is a linear, self-adjoint elliptic operator, and therefore has a discrete spectrum of eigenvalues, infinite in both directions. For $\nu \gg 0$, let us denote by $V^{\nu}$ the direct sum of all eigenspaces of $\ell$ with eigenvalues in the interval $(-\nu, \nu]$. The spaces $V^{\nu}$ are finite dimensional and preserved by the map $\ell$. As $\nu \to \infty$, these spaces provide a finite dimensional approximation for $V$, in the sense that the $L^2$-projection $p^{\nu}: V \to V^{\nu} \subset V$ limits to the identity $\id_V$ pointwise. 

Instead of considering the flow trajectories of $\ell + c$ on $V$, we will look at trajectories of $\ell + p^{\nu} c$ on $V^{\nu}$. (These are gradient flow trajectories for the restriction of $\csd$ to $V^{\nu}$, in a suitable metric.) As discussed in Section~\ref{sec:Conley}, to be able to apply Conley index theory we need to make sure that the critical points and flow lines between them form a compact set. This is true for the original Seiberg-Witten equations on $V$, by standard results in gauge theory. In order to obtain the same result for the approximate flow $\varphi^{\nu}$, we must restrict to an a priori bounded set (in an $L^2_k$ norm); this is because the projections $p^{\nu}$ do not converge to $1$ strongly in a Sobolev norm as $\nu \to \infty$; they only do so pointwise (and hence uniformly on compact sets). Let $R$ be the a priori $L^2_k$ bound on the size of Seiberg-Witten solutions on $V$. If we restrict to a larger ball, say $B(2R)$, then it can be shown that $\ell + p^{\nu} c$ converges to $\ell + c$ uniformly there, and therefore, for large $\nu$, the solutions to $\ell + p^{\nu}c = 0$ are inside the smaller  ball $B(R)$. A similar argument applies to points on the flow lines connecting critical points. We deduce that if we define $\S^{\nu}$ as the set of all critical points and flow lines of $\ell + p^{\nu}c$ in $B(2R) \cap V^{\nu}$, then $\S^{\nu}$ is compact. Further, the flow $\ell + p^{\nu}c$ and the set $S^{\nu}$ are $S^1$-invariant. Therefore, as explained in Section~\ref{sec:Conley}, there is an associated $S^1$-equivariant Conley index $I^{\nu} := I_{S^1}(\varphi^{\nu}, \S^{\nu})$. 
  
We define the $S^1$-equivariant Seiberg-Witten Floer homology of $(Y, \s)$ to be the (reduced) equivariant homology of $I^{\nu}$, with a shift in degree:
\begin{equation}
\label{eq:swfh}
 \swfh^{S^1}_*(Y, \s) := \tH_{*+ \dim_{\R} V^0_{-\nu} +2n(Y, \s, g)}^{S^1}(I^\nu).
 \end{equation}

Here, $V^0_{-\nu}$ stands for the direct sum of the eigenspaces of $\ell$ with eigenvalues between $-\nu$ and $0$, and $n(Y, \s, g) \in  \Q$ is a certain quantity (a combination of eta invariants) depending on the metric $g$ on $Y$. The shift by $\dim V^0_{-\nu}$ is necessary because as we change $\nu$ to some $\nu' > \nu$, the Conley index changes by a suspension: $I^{\nu'} = (V_{-\nu'}^{-\nu})^+ \wedge I^{\nu}$. The shift by $n(Y, \s, g)$ is needed to compensate for the change in the dimension of $V^0_{-\nu}$ as we vary the metric $g$. If we have a family of metrics $(g_t)_{t \in [0,1]}$, then $n(Y, \s, g_0) - n(Y, \s, g_1)$ is the spectral flow of $\ell$ in that family. (In fact, it is the spectral flow of $\dirac$, because we assumed $b_1(Y)=0$ and hence $*d$ has trivial spectral flow.)

We mention that $n(Y, \s, g)$ can be computed as follows. Choose a compact $4$-manifold $W$ with boundary $Y$, and let $\t$ be a $\spinc$ structure on $W$ that restricts to $\s$ on $Y$. Equip $W$ with a Riemannian metric such that a neighborhood of the boundary is isometric to $[0,1] \times Y$, and let $\Dirac$ be the Dirac operator on $(W, \t)$ with spectral boundary conditions as in \cite{APS}. Then:
\begin{equation}
\label{eq:n}
 n(Y, \s, g) =  \ind_{\C}(\Dirac) - (c_1(\t)^2 - \sigma(W))/8.
 \end{equation}

\begin {remark}
The degree shift in \eqref{eq:swfh} has a parallel in the versions of Seiberg-Witten Floer homology defined Morse-theoretically. If we have a complex $\swfc_*^{S^1}(Y, \s, g)$ as in Section~\ref{sec:sw}, 
then one defines an absolute grading on it by setting the lowest group in the $U$-tail \eqref{eq:Utail} to be in degree $-2n(Y, \s, g)$.  
\end{remark}

Using finite dimensional approximation we can define a more refined invariant than Floer homology. Recall that the Conley index $I^{\nu}$ is a homotopy type. When we vary $\nu$ this changes by suspensions. We can introduce formal de-suspensions of $I^{\nu}$ to produce an invariant of $(Y, \s)$ in the form of an $S^1$-equivariant based stable homotopy type:
\begin{equation}
\label{eq:swf}
\swf(Y, \s) := \Sigma^{-V^0_{-\nu}} \Sigma^{-n(Y, \s, g) \C} I^\nu.
\end{equation}
Here, $\C$ denotes a copy of the standard one-dimensional complex representation of $S^1$. The $S^1$-equivariant homology of $\swf(Y, \s)$ is the Seiberg-Witten Floer homology defined in \eqref{eq:swfh}.

Let $\Cat$ be the $S^1$-equivariant analog of the Spanier-Whitehead category of suspension spectra. By keeping careful track of the orientations of eigenspaces of $\ell$, one can define $\swf(Y, \s)$ 
as an element of $\Cat$, up to canonical equivalence; compare Remark~\ref{rem:css}.

\subsection{A Pin(2)-equivariant version} \label{sec:pin}
The group $\pin$ (sometimes known as $\pin^-$) is a non-trivial extension of $\Z/2$ by $S^1$. An easy way to define it is as a subgroup of the unit quaternions $S(\H) \cong \sutwo$: If we write $\H$ as $\C \oplus \C j$, then $\pin=S^1 \cup S^1 j$. 

The Seiberg-Witten equations are invariant under conjugation of $\spinc$ structures: $\s \mapsto \bar \s$. If we combine this with the $S^1$-action, we obtain a $\pin$-action. This is particularly interesting when $\s$ comes from a spin structure, so that $\s = \bar \s$. The spinor bundle $\Spin$ is then quaternionic, and the action of $j \in \pin \subset S(\H)$ on $(a, \phi)$ can be written as 
$$j : (a, \phi) \mapsto(-a, \phi j).$$ 

Incorporating the $\pin$-symmetry into a Morse-theoretic approach to Seiberg-Witten Floer homology seems rather difficult. On the other hand, doing it within the context of finite dimensional approximation is straightforward. The Conley index $I^{\nu}$ can be taken to be $\pin$-equivariant, and we define the 
$\pin$-equivariant Seiberg-Witten Floer homology of $(Y, \s)$ as
\begin{equation}
\label{eq:swfhpin2}
 \swfh^{\pin}_*(Y, \s) := \tH_{*+ \dim V^0_{-\nu} +2n(Y, \s, g)}^{\pin}(I^\nu).
 \end{equation}
We can also define a $\pin$ version of $\swf(Y, \s)$, as the $\pin$-equivariant stable homotopy type of $I^{\nu}$ with the same formal de-suspension as before.

Observe that $\swfh^{\pin}_*(Y, \s)$ is a module over the $\pin$-equivariant cohomology of a point, i.e. $H^*(B \pin)$. To compute $H^*(B \pin)$, we use the fiber bundle
  $$\pin \To \sutwo \To \rp^2$$
which yields another fiber bundle relating the classifying spaces:
$$ \rp^2 \To B\pin \To B\sutwo = \hp^{\infty}.$$
The associated Leray-Serre spectral sequence has no room for higher differentials. Consequently, $H^*(B \pin)$
is isomorphic to $H^*(\hp^{\infty}) \otimes H^*(\rp^2)$. 

For simplicity, let us work with coefficients in the field $\F$ with two elements. Then $H^*(B\pin; \F) = \F[q,v]/(q^3)$, with $q$ in degree $1$ and $v$ in degree $4$. It is helpful to imagine that $\swfh^{\pin}_*(Y, \s; \F)$ is the homology of a complex $\swfc^{\pin}_*(Y, \s; \F)$, composed of a copy of $H^*(B\pin; \F)$ for the reducible Seiberg-Witten solution, and a copy of $\F$ for each pair of ($S^1$-orbits of) irreducible solutions related by $j$. Thus, the reducible contributes a triple of infinite $v$-tails:
\begin{equation}
\label{eq:vtails}
 \xymatrixcolsep{.7pc}
\xymatrix{
 \F  &  \F \ar@/_1pc/[l]_{q} &  \F \ar@/_1pc/[l]_{q} & 0 & \F \ar@/^1pc/[llll]^{v} & \F \ar@/_1pc/[l]_{q} \ar@/^1pc/[llll]^{v} & \F \ar@/_1pc/[l]_{q} \ar@/^1pc/[llll]^{v} & 0 & \dots  \ar@/^1pc/[llll]^{v} & \dots \ar@/^1pc/[llll]^{v} & \dots \ar@/^1pc/[llll]^{v}
} 
\end{equation}
The generators in the complex are related to each other by $\del$, $v$ and $q$ maps, lowering degrees by $1$, $4$ and $1$, respectively. The absolute grading is again obtained by requiring the lowest degree generator in \eqref{eq:vtails} to be in grading $-2n(Y, \s, g)$.

The existence of the complex $\swfc^{\pin}_*(Y, \s; \F)$ is at this point only a heuristic, since we don't have a Morse-theoretic definition of $\pin$-equivariant Seiberg-Witten Floer homology. Nevertheless, in some cases (such as for Brieskorn spheres) one can pick a metric and perturbation such that there is a complex of that form, and show that its homology is $\swfh^{\pin}_*(Y, \s; \F)$.

\subsection{Examples} \label{ex:ex} The following three examples are integral homology spheres. In such cases there is a unique $\spinc$ structure $\s$, which we drop from the notation. Further, for integral homology spheres, the quantity $n(Y, g)$ is always an integer.

The simplest example is $S^3$. If we pick $g$ to be the round metric, then there are no irreducibles, and the value of $n(S^3, g)$ is zero. Thus, the $S^1$- and $\pin$-equivariant Seiberg-Witten Floer homologies of $S^3$ look exactly like \eqref{eq:Utail} and \eqref{eq:vtails}, respectively, with the lowest terms in degree $0$. The $\pin$-equivariant stable homotopy type $\swf(S^3)$ is that of $S^0$, with trivial $\pin$ action.

The case of the Poincar\'e sphere $P=\Sigma(2,3,5)$ (with its round metric) is very similar, with no irreducibles, but with the difference that we have $n(P, g)=-1$. Thus, the lowest term in the two Floer homologies is in degree $2$.

For a more non-trivial example, consider the Brieskorn sphere $\Sigma(2,3,11)$. Using a metric as in \cite{MOY}, we find that the irreducibles form two $S^1$-orbits, related by the action of $j$ ; in other words, they form one $\pin$-orbit. The value of $n(\Sigma(2,3,11), g)$ is still zero. The irreducibles are in degree $1$ and they interact with the reducible through the $\del$ map. Thus, the $S^1$-equivariant Seiberg-Witten Floer complex of $\Sigma(2,3,11)$ is
\begin{equation}
\label{eq:11c}
 \xymatrixcolsep{.7pc}
\xymatrixrowsep{0pc}
\xymatrix{
  \Z  &  0 & \Z \ar@/_1pc/[ll]_{U} & 0 & \Z \ar@/_1pc/[ll]_{U} & 0 & \dots \ar@/_1pc/[ll]_{U} \\
& \oplus  & & & & & &  \\
  \ \ \  & \Z\ar[uul]^{\del} & & & & & & \\
 & \oplus & & & & & &  \\
\ \ \  & \Z\ar@/^1pc/[uuuul]^{\del} & & & & & &  
} \end{equation}
with the leftmost element in degree $0$. Its homology is
\begin{equation}
\label{eq:11h}
 \xymatrixcolsep{.7pc}
\xymatrixrowsep{0pc}
\xymatrix{
{\phantom{\Z}} &  0 & \Z  & 0 & \Z \ar@/_1pc/[ll]_{U} & 0 & \dots \ar@/_1pc/[ll]_{U} \\
& \oplus  & & & & & &  \\
  \ \ \  & \Z & & & & & & \\
 } \end{equation}
with the bottom $0 \oplus \Z$ in degree $1$.
 
The $\pin$-equivariant Seiberg-Witten Floer complex of $\Sigma(2,3,11)$ is
\begin{equation}
\label{eq:pin11c}
 \xymatrixcolsep{.7pc}
\xymatrixrowsep{0pc}
\xymatrix{
  \F  &  \F \ar@/_1pc/[l]_{q} &  \F \ar@/_1pc/[l]_{q} & 0 & \F \ar@/^1pc/[llll]^{v} & \F \ar@/_1pc/[l]_{q} \ar@/^1pc/[llll]^{v} & \F \ar@/_1pc/[l]_{q} \ar@/^1pc/[llll]^{v} & 0 & \dots  \ar@/^1pc/[llll]^{v} & \dots \ar@/^1pc/[llll]^{v} & \dots \ar@/^1pc/[llll]^{v}\\
   & \oplus  & & & & & & & & & \\
    & \F\ar@/^1pc/[uul]^{\del} & & & & & & & & & 
} \end{equation}
again with the leftmost element in degree $0$. The homology
\begin{equation}
\label{eq:pin11h}
 \xymatrixcolsep{.7pc}
\xymatrixrowsep{0pc}
\xymatrix{
{\phantom{\F}} &  \F  &  \F \ar@/_1pc/[l]_{q} & 0 & \F  & \F \ar@/_1pc/[l]_{q} \ar@/^1pc/[llll]^{v} & \F \ar@/_1pc/[l]_{q} \ar@/^1pc/[llll]^{v} & 0 & \dots  \ar@/^1pc/[llll]^{v} & \dots \ar@/^1pc/[llll]^{v} & \dots \ar@/^1pc/[llll]^{v} 
} \end{equation}
has the leftmost element in degree $1$.

The stable homotopy type $\swf(\Sigma(2,3,11))$ is that of the unreduced suspension of $\pin$, with one of the cone points as the basepoint, and with the induced $\pin$-action.

\subsection{Properties}
Let us now describe a few properties of Seiberg-Witten Floer homologies and stable homotopy types. We will omit the $\spinc$ structures from notation for simplicity. 

\subsubsection{Orientation reversal} If we change the orientation of $Y$, then the approximate Seiberg-Witten flow $\varphi^{\nu}$ changes direction. Let us recall from Section~\ref{sec:ci} that in this case the two Conley indices are Spanier-Whitehead dual. (This is also true equivariantly.) This implies that $X=\swf(Y)$ and $X'=\swf(-Y)$ are dual in the equivariant Spanier-Whitehead category.

Non-equivariantly, if $X$ and $X'$ are dual to each other, then the homology of $X$ is isomorphic to the cohomology of $X'$, with the degrees changing sign: $\tH_*(X) \cong \tH^{-*}(X')$. Equivariantly, this cannot be true exactly as such, because $\tH_*(X)$ is unbounded in the positive direction (with regard to grading), but bounded below, whereas $\tH^{-*}(X')$ is bounded above but not below. Instead, what happens is that (for any group $G$, in our case $S^1$ or $\pin$) we have a long exact sequence
\begin{equation}
\label{eq:tate}
 \dots \To \tH^{-*}_G(X') \To t\tH_{*}^G(X) \To \tH_{*-\dim G - 1}^G(X) \To  \dots 
\end{equation}
Here, $t\tH_{*}^G(X)$ denotes the $G$-equivariant Tate homology of $X$, which depends only mildly on $X$. (When $G$ is trivial, the Tate homology is zero.) In our setting, the $S^1$-equivariant Tate homology of $\swf(Y)$ is always isomorphic to $\Z[U, U^{-1}]$, that is, to an infinite $U$-tail in both directions:
\begin{equation}
\label{eq:tate2}
 \xymatrixcolsep{.7pc}
\xymatrix{
 \dots & 0 & \Z \ar@/_1pc/[ll]_{U}  &  0 & \Z \ar@/_1pc/[ll]_{U} & 0 & \Z \ar@/_1pc/[ll]_{U} & 0 & \dots \ar@/_1pc/[ll]_{U} 
} \end{equation}

Similarly, the $\pin$-equivariant Tate homology of $\swf(Y)$ is always isomorphic to $\F[q, v, v^{-1}]/(q^3).$

\subsubsection{Disjoint unions} \label{sec:du} The Seiberg-Witten Floer stable homotopy type of a disjoint union $Y_0 \amalg Y_1$ is the smash product
$$ \swf(Y_0 \amalg Y_1) = \swf(Y_0) \wedge \swf(Y_1).$$
Hence, the corresponding Seiberg-Witten Floer cohomologies are related by K\"unneth spectral sequences.

\subsubsection{Cobordism} \label{sec:cobordism} A very important feature of Floer homology is that it fits into a form of TQFT (topological quantum field theory). Let $G$ be $S^1$ or $\pin$. Suppose we have smooth four-dimensional cobordism $W$ from $Y_0$ to $Y_1$, equipped with a $\spinc$ structure $\t$ in the case $G=S^1$, or a spin structure $\t$ in the case $G=\pin$. Then, the Seiberg-Witten equations on $W$ produce a module homomorphism from $\swfh_*^G(Y_0)$ to $\swfh_*^G(Y_1)$, with a shift in degree depending on $W$ and $\t$. Moreover, these homomorphisms come from an actual map of suspension spectra $\swf(Y_0) \to \swf(Y_1)$. The maps are functorial with respect to composition of cobordisms $Y_0 \xrightarrow{W_0} Y_1 \xrightarrow{W_1} Y_2$, as long as the middle manifold $Y_1$ is connected. We refer to \cite{GluingBF} for more details.

\section{The Triangulation Conjecture}
\label{sec:tc}
\subsection{Background}
A natural question in topology is whether every manifold admits a simplicial triangulation, that is, a homeomorphism to a simplicial complex. A triangulation would allow the manifold to be described in simple combinatorial terms. The origins of the triangulation problem go back to Poincar\'e, who gave an incomplete proof in Chapter XI of his first supplement to Analysis Situs \cite{Poincare, PoincareBook}. Poincar\'e was working with differentiable manifolds (although the terminology and the rigor were not there yet). Much later, Cairns \cite{Cairns} and Whitehead \cite{Whitehead} showed that every differentiable manifold can indeed be triangulated.

In 1924, Kneser \cite{Kneser} asked the triangulation question for topological manifolds. The answer was thought to be positive, and this became known as the Triangulation Conjecture. However, in the end the conjecture turned out to be false.  More precisely, the answer depends on the dimension of the manifold:
\begin{itemize}
\item The conjecture is true in dimension $\leq 3$. This was proved for surfaces by Rad\'o \cite{Rado}, and for three-manifolds by Moise \cite{Moise}. (Dimensions zero and one are easy.)
\item It is false in dimension $4$. This is where the first counterexamples were found. In the mid 1980's, Casson \cite{Casson} showed that, for example, Freedman's $E_8$-manifold \cite{Freedman} is non-triangulable.
\item It is also false in dimensions $\geq 5$. The disproof consists of two parts. The first, due to Galewski-Stern \cite{GS} and, independently, Matumoto \cite{Matumoto}, reduces the problem to a question about cobordisms of three-manifolds. The second part is the solution to this question \cite{beta}, which uses $\pin$-equivariant Seiberg-Witten Floer homology and will be sketched shortly.
\end{itemize}

Let us also mention here a related, stronger question: Does every topological manifold admit a combinatorial triangulation, that is, a piecewise linear (PL) structure? A triangulation is called  {\em combinatorial} if the links of the vertices are PL-homeomorphic to spheres. The answers (depending on dimension) in the combinatorial case are the same as before, but the chronology of discoveries was different:
\begin{itemize}
\item Every manifold of dimension $\leq 3$ has a PL structure. (The triangulations found by Rad\'o and Moise were combinatorial.)
\item Kirby and Siebenmann \cite{KSbook} showed that there exist manifolds without PL structures in every dimension $\geq 5$. Specifically, they showed that every topological manifold $M$ has an associated obstruction class $\Delta(M) \in H^4(M; \Z/2)$, and $M$ has a PL structure if and only if $\Delta(M)=0$. Further, they showed that for any $n \geq 5$, there exist $n$-dimensional manifolds $M$ with $\Delta(M) \neq 0$.
\item Freedman \cite{Freedman} found $4$-dimensional manifolds without PL structures. An example is his $E_8$-manifold: a closed, simply connected $4$-manifold with intersection form $E_8$. Indeed, in  dimension four PL structures are equivalent to smooth ones. A smooth, simply connected $4$-manifold with even intersection form is spin, and hence has signature divisible by $16$ by Rokhlin's theorem \cite{Rokhlin}; therefore, the $E_8$-manifold is not PL. 
\end{itemize}

We point out that for closed, oriented, spin $4$-manifolds the Kirby-Siebenmann class is given by 
$$ \Delta(M) = \sigma(M)/8 \! \pmod 2 \ \in \Z/2 \cong H^4(M; \Z/2). $$
The Kirby-Siebenmann class still obstructs PL structures in dimension four; for example, the $E_8$-manifold has $\Delta \neq 0$. However, in this dimension there are also other obstructions, coming from gauge theory; e.g. Donaldson's diagonalizability theorem \cite{Donaldson} and Furuta's $10/8$ theorem \cite{Furuta}. For instance, the connected sum of two copies of the $E_8$ manifold has $\Delta = 0$, but is non-smoothable (hence not PL) by Donaldson's theorem.

In dimensions $n \geq 5$, an example of a non-PL manifold is the product of the $E_8$-manifold with the torus $T^{n-4}$. Of course, the original examples of Kirby and Siebenmann were different, since they came before Freedman's work; see \cite{KSbook} for their construction.

We refer to Ranicki's survey \cite{Ranicki} for more details about this subject, and about the related  Hauptvermutung.

It is worth mentioning that every topological manifold $M$ is homotopy equivalent to a simplicial complex.\footnote{Laurence Taylor informed the author that, in fact, the simplicial complex can be taken to be of the same dimension as the manifold. In the published version of this paper, this fact was mistakenly listed as an open problem.} (Moreover, if the manifold is compact, then the simplicial complex can be taken to be finite.) This is a consequence of the work of Kirby and Siebenmann \cite{KSbook}. One starts by embedding the manifold into a large Euclidean space $\R^m$. The associated sphere bundle (the boundary of a standard neighborhood of $M$) has a trivial normal line bundle, and this implies that it can be isotoped to a PL submanifold of $\R^m$. It follows that the associated disk bundle (which is homotopy equivalent to $M$) admits a PL structure. 

Lastly, if we weaken ``simplicial complex'' to ``CW complex,'' then it is known that every topological manifold of dimension $d \neq 4$ has a handlebody structure, and hence is homeomorphic to a CW complex. See \cite[p. 104]{KSbook} for the case $d > 5$ and \cite{Quinn} for the case $d=5$. It is an open problem whether every manifold of dimension $4$ is homeomorphic to a CW complex.

\subsection{Reduction to a question about homology cobordism}
We now return to the question of the existence of non-triangulable manifolds in dimensions $\geq 5$. Note that most triangulations of manifolds that one can think of are combinatorial. In fact, for a long time it was not known if any non-combinatorial triangulations existed. This changed with the work of Edwards \cite{Edwards}, who provided the first examples.

\begin{example}
Let $K$ be a triangulation of a non-trivial homology sphere $M$ of dimension $n \geq 3$, such that $\pi_1(M) \neq 1$; for example, $M$ could be the Poincar\'e sphere. Consider the suspension $\Sigma M$ of $M$; the triangulation $K$ induces one on $\Sigma M$. The space $\Sigma M$ is not a manifold, because if we delete a cone point $x$ from a neighborhood of $x$, then the result is not simply connected. However, if we repeat the procedure and construct the double suspension $\Sigma^2 M$, the Double Suspension Theorem of Edwards \cite{Edwards, EdwardsICM} and Cannon \cite{Cannon} tells us that $\Sigma^2 M$ is homeomorphic to $S^{n+2}$. The induced triangulation on $\Sigma^2 M \cong S^{n+2}$ is not combinatorial, because the links of the two cone points are not spheres. 
\end{example}

\begin{remark}
In dimensions $n \leq 4$, every simplicial triangulation of an $n$-manifold is combinatorial: The link of every vertex can be shown to be a simply connected, closed $(n-1)$-manifold, and therefore to be the $(n-1)$-sphere. (For $n=4$, this argument uses the Poincar\'e Conjecture, proved by Perelman \cite{Perelman1, Perelman2, Perelman3}.) 
\end{remark}

Let us now suppose that a closed, oriented $n$-dimensional manifold $M$ ($n \geq 5$) is equipped with a triangulation $K$. Consider the following element (sometimes called the Sullivan-Cohen-Sato class; cf. \cite{Sullivan, Cohen, Sato}):
\begin{equation}
\label{eq:cK}
 c(K) = \sum_{\sigma \in K^{(n-4)}} [\link_K(\sigma)] \cdot \sigma \in H_{n-4}(M; \Theta^H_3) \cong H^4(M; \Theta^H_3).
 \end{equation}
Here, the sum is taken over all codimension four simplices in the triangulation $K$. The link of each such simplex can be shown to be a homology $3$-sphere. (It would be an actual $3$-sphere if the triangulation were combinatorial.) The group $\Theta^H_3$ is the three-dimensional homology cobordism group; it is generated by equivalence classes of oriented integral homology $3$-spheres,  where $Y_0$ is equivalent to $Y_1$ if there exists a piecewise-linear (or, equivalently, a smooth) compact, oriented $4$-dimensional cobordism $W$ from $Y_0$ to $Y_1$, such that $H_1(W; \Z) = H_2(W; \Z)=0$. Addition in $\Theta^H_3$ is given by connected sum, the inverse is given by reversing the orientation, and $S^3$ is the zero element.  

The reader may wonder why we focus on codimension four simplices in \eqref{eq:cK}. The reason is that the analog of the homology cobordism group in any other dimension is trivial: For $n \neq 3$, every $n$-dimensional PL homology sphere is the boundary of a contractible PL manifold, according to a theorem of Kervaire \cite{Kervaire}.

In dimension three, it is known that $\Theta^H_3$ is infinite, and in fact infinitely generated \cite{FSorbifolds, FurutaHom, FSinstanton}. However, its general structure is still a mystery: for example, it is not known whether $\Theta^H_3$ has any non-trivial torsion elements. 

In the study of triangulations, an important role is played by the Rokhlin homomorphism \cite{Rokhlin, EellsKuiper}:
$$\mu: \Theta^H_3 \to \Z/2, \ \  \mu(Y) = {\sigma(W)}/{8} \pmod {2},$$
where $W$ is any compact, spin $4$-manifold with boundary $Y$. Rokhlin's theorem shows that the value of $\mu$ depends only on $Y$, not on $W$. The homomorphism $\mu$ can be used to show that $\Theta^H_3$ is non-trivial: For instance, the Poincar\'e sphere $P$ bounds the $E_8$ plumbing (of signature $-8$), so $\mu(P)=1$.

Consider the short exact sequence
\begin{equation}
\label{eq:ses}
 0 \To \ker(\mu) \To \Theta^H_3 \To \Z/2 \To 0
 \end{equation}
and the associated long exact sequence in cohomology
\begin{equation}
\dots \To H^4(M; \Theta^H_3) \xrightarrow{\hspace*{2pt} \mu_* \hspace*{2pt}} H^4(M; \Z/2)  \xrightarrow{\hspace*{2pt} \delta \hspace*{2pt}} H^5(M; \ker(\mu)) \To \dots,
\end{equation}
where $\delta$ denotes the Bockstein homomorphism.

Let us return to the element $c(K)$ defined in \eqref{eq:cK}. Clearly, if the triangulation $K$ is PL, then $c(K) = 0$. It can be shown that the image of $c(K)$ under $\mu_*$ is exactly the Kirby-Siebenmann obstruction to PL structures, $\Delta(M) \in H^4(M; \Z/2)$. (This gives a simple way of thinking about $\Delta$, albeit one that only applies to triangulable manifolds.) Thus, $M$ admits a PL triangulation (possibly different from $K$) if and only if $\mu_*(c(K)) = \Delta(M) =0$. 

If $M$ admits any triangulation, we get that $\Delta(M)$ is in the image of $\mu_*$, and hence in the kernel of $\delta$. Thus, a necessary condition for the existence of simplicial triangulations is the vanishing of the class
$$ \delta(\Delta(M)) \in  H^5(M; \ker(\mu)).$$

It can be shown that this is also a sufficient condition. Further, while the discussion above (inspired from   \cite{Ranicki}) was for the case when $M$ closed and oriented, these assumptions are not necessary for the conclusion:

\begin{theorem}[Galewski-Stern \cite{GS}; Matumoto \cite{Matumoto}]
A topological manifold $M$ of dimension $\geq 5$ is triangulable if and only if $\delta(\Delta(M)) =0$.
\end{theorem}

We are left with the question of whether there exist $M$ with $\delta(\Delta(M))\neq 0$. Observe that the Bockstein homomorphism $\delta$ is guaranteed to vanish if the short exact sequence \eqref{eq:ses} splits. Thus, if \eqref{eq:ses} splits, then all high dimensional manifolds would be triangulable. In fact, we have:

\begin{theorem}[Galewski-Stern \cite{GS}; Matumoto \cite{Matumoto}]
\label{thm:gsm}
There exist non-triangulable manifolds of (every) dimension $\geq 5$ if and only if the exact sequence \eqref{eq:ses} does not split.
\end{theorem}

A few remarks are in order about the ``if'' part of the theorem. Galewski and Stern \cite{GS5} constructed an explicit five-dimensional manifold $M$ with $\Sq^1 \Delta(M) \neq 0 \in H^5(M; \Z/2)$, where $\Sq^1$ denotes the first Steenrod square. The first Steenrod square is the Bockstein homomorphism for the exact sequence
$$ 0 \To \Z/2 \To \Z/4 \To \Z/2 \To 0$$
and a little algebra shows that if \eqref{eq:ses} does not split, then the non-vanishing of $\Sq^1 \Delta(M)$ implies the non-vanishing of $\delta \Delta(M) \in H^5(M; \ker(\mu))$.

Moreover, if $M$ is a five-manifold with $\Sq^1 \Delta(M) \neq 0$, then the products $M \times T^{n-5}$ provide examples of manifolds with the same property in every dimension $\geq 5$.

It turns out that \eqref{eq:ses} does not split (cf. Theorem~\ref{thm:beta} below), so the Galewski-Stern manifold from \cite{GS5} is non-triangulable. Here is a different example, based on Freedman's work on four-manifolds:

\begin{example}[Peter Kronheimer]
By Freedman's theorem \cite{Freedman}, simply connected, closed topological four-manifolds are characterized (up to homeomorphism) by their intersection form and their Kirby-Siebenmann invariant. Let $W$ be the fake $\cp^2 \# (-\cp^2)$, that is, the closed, simply connected topological $4$-manifold with intersection form $Q=\langle 1 \rangle \oplus \langle -1 \rangle$ and non-trivial Kirby-Siebenmann invariant. Since the form $Q$  is isomorphic to $-Q$, by applying Freedman's theorem again we find that $W$ admits an orientation-reversing homeomorphism $f: W \to W$. Let $M$ be the mapping torus of $f$. Then $M$ is a five-manifold with
$w_1(TM)$ Poincar\'e dual to the class $[W \times \pt]$. Further, $\Delta(M) \in H^4(M; \Z/2)$ is Poincar\'e dual to the class of a section of the bundle $M \to S^1$. Therefore, by Wu's formula,
$$ \Sq^1 \Delta(M) = \Delta(M) \cup w_1(TM) =1 \in H^5(M; \Z/2) \cong \Z/2.$$
We deduce that $M$ is non-triangulable.
\end{example}

\begin{example}
If instead of the mapping torus of $f$ we would simply consider the manifold $M' = W \times S^1$, then we would get $\Delta(M')\neq 0$ but $\delta(\Delta(M'))=0$. Thus, $M'$ admits simplicial triangulations but not combinatorial triangulations. 
\end{example}

By the work of Siebenmann \cite[Theorem B]{SiebenmannTC} combined with the Double Suspension Theorem \cite{Edwards, Cannon}, it follows that all $5$-dimensional non-triangulable manifolds have to be compact and non-orientable. However, this is not the case in higher dimensions:

\begin{example}[Ron Stern]
Let $M$ be a non-triangulable five-dimensional manifold.  Since $M$ is necessarily non-orientable, it admits an oriented double cover $\tM \to M$. Consider the six-dimensional manifold
$$ N = \tM \times_{\Z/2} S^1$$
which is an unoriented $S^1$-bundle over $M$. Then $N$ is orientable, but non-triangulable because we still have $\delta (\Delta(N)) \neq 0$. We can also get a non-compact example by replacing $S^1$ with $\R$ in this construction.
\end{example}

In another direction, Davis, Fowler and Lafont have shown that in every dimension $\geq 6$ there exist non-triangulable aspherical manifolds \cite{DFL}.

\subsection{Solution using Seiberg-Witten theory}
\label{sec:beta}
Theorem~\ref{thm:gsm} reduced the triangulation problem in high dimensions to a question about the group $\Theta^H_3$. In this section we sketch its solution:

\begin{theorem}[\cite{beta}]
\label{thm:beta}
The short exact sequence \eqref{eq:ses} does not split.
\end{theorem} 

A splitting of \eqref{eq:ses} would consist of a map $\eta: \Z/2 \to \Theta^H_3$ with $\mu\circ \eta=\id$; that is, we would need a homology $3$-sphere $Y$ such that $Y$ has Rokhlin invariant one, and $Y$ is of order two in the homology cobordism group. 

To show that such a sphere does not exist, it suffices to construct a lift of $\mu$ to the integers,
$$ \beta : \Theta^H_3 \to \Z,$$
with the following properties:
\begin{enumerate}[(a)]
\item If $-Y$ denotes $Y$ with the orientation reversed, then $\beta(-Y) = - \beta(Y)$;
\item The mod $2$ reduction of $\beta(Y)$ is the Rokhlin invariant $\mu(Y)$.
\end{enumerate}

We will construct a map $\beta$ of this type. Interestingly, this map will not be a homomorphism. (For more on this, see Section~\ref{sec:NoAdd} below.) Nevertheless, properties (a) and (b) suffice to prove Theorem~\ref{thm:beta}. Indeed, if we had a  homology sphere $Y$ of order two in $\Theta^H_3$, then $Y$ would be homology cobordant to $-Y$, and we would obtain
$$ \beta(Y) = \beta(-Y) = - \beta(Y),$$
hence $\beta(Y)=0$ and therefore $\mu(Y) = 0$.

The construction of $\beta$ involves $\pin$-equivariant Seiberg-Witten Floer homology. The definition can be extended to rational homology spheres equipped with spin structures; in that case $\beta$ takes values in $\tfrac{1}{8} \Z \subset \Q$, rather than in $\Z$. For simplicity, we will only discuss the case of integral homology spheres.

Before explaining $\beta$, let us recall a predecessor, the Fr{\o}yshov invariant from \cite{FroyshovSW, KMOS, KMbook}. (A parallel construction exists in Heegaard Floer homology \cite{AbsGraded}.) The Fr{\o}yshov invariant is defined from $S^1$-equivariant Seiberg-Witten Floer homology. Suppose that $Y$ is an integral homology sphere. Recall that $\swfh^{S^1}_*(Y)$ is the direct sum of an infinite $U$-tail as in \eqref{eq:Utail} and a finitely generated piece. The Fr{\o}yshov invariant $h(Y)$ is defined as $-d(Y)/2$, where $d(Y)$ is the minimal grading of an element in the $U$-tail. (So that there is no confusion about whether we allow torsion elements in the tail, it is convenient to fix a field $\f$ and work with coefficients in $\f$ rather than $\Z$.) It is important to note that we consider the $U$-tail in Floer homology, not in the chain complex $\swfc^{S^1}_*(Y, g)$. For the tail in the chain complex, the minimal degree is the quantity $-2n(Y, g)$, which depends on the metric $g$. In  contrast, $d(Y)$ is an invariant of $Y$.

The Fr{\o}yshov invariant descends to a map
$$ h: \Theta^H_3 \to \Z.$$
The proof of this fact is based on the TQFT properties of $\swfh^{S^1}$ discussed in Section~\ref{sec:cobordism}. Precisely, if $Z$ is a homology cobordism from $Y_0$ to $Y_1$, then there is an induced map $\swfh^{S^1}_*(Y_0) \to \swfh^{S^1}_*(Y_1)$, without a shift in degree. Moreover, the map is an isomorphism between the infinite $U$-tails in large enough degrees. (This can be seen by studying reducible Seiberg-Witten solutions on $Z$.) The module structure of $\swfh^{S^1}_*$ implies that the bottom degree of the $U$-tail for $Y_0$ needs to be smaller or equal to the bottom degree of the $U$-tail for $Y_1$. By reversing the cobordism, we get an inequality in the opposite direction, so we can conclude that $h(Y_0) = h(Y_1)$.
 
The Fr{\o}yshov invariant satisfies the analog of property (a) for $\beta$, that is, we have $h(-Y) = -h(Y)$. This can be proved using the long exact sequence \eqref{eq:tate}. Given that this is an exact sequence of $\Z[U]$-modules, we see that the Tate homology \eqref{eq:tate2} must be composed of the $U$-tail for $Y$ together with the reverse of the one for $-Y$. This implies that $h(-Y) = -h(Y)$.

However, the Fr{\o}yshov invariant does not reduce mod $2$ to the Rokhlin invariant. This is in spite of the fact that the minimal degree $-2n(Y, g)$ of the $U$-tail on the chain complex does capture the Rokhlin invariant. Indeed, recall from \eqref{eq:n} that we have $n(Y, g) = \ind_{\C}(\Dirac)- (c_1(\t)^2-\sigma(W))/8$. Since $H_1(Y)=0$, we can choose $\t$ to be a spin structure, so that $c_1(\t)=0$ and the Dirac operator $\Dirac$ acts on a quaternionic vector space. We get that $\ind_{\C}(\Dirac)=2\ind_{\H}(\Dirac)$ is even and hence $n(Y, g) = \ind_{\C}(\Dirac)+ (\sigma(W)/8)$ reduces to $\mu(Y)$ modulo $2$. Therefore, the minimal degree of the $U$-tail on the chain complex is an even integer, congruent to $2\mu(Y)$ modulo $4$. The minimal degree of the $U$-tail on homology, $d(Y)$, is still an even integer, but because of the interaction with the irreducibles we cannot say much about its congruence class mod $4$. Hence, the parity of $h(Y)=-d(Y)/2$ is unrelated to $\mu(Y)$.

\begin{example}
Consider the Brieskorn sphere $\Sigma(2,3,11)$. The bottom $\Z$ in the $U$-tail from the chain complex \eqref{eq:11c} is in degree zero, in agreement with the fact that $\mu(\Sigma(2,3,11))=0$. The differential $\del$ cancels this $\Z$ against one $\Z$ from the irreducibles. Thus, the bottom $\Z$ in the $U$-tail in the homology \eqref{eq:11h} is in degree $2$. We get that $d(\Sigma(2,3,11))=2$, and that $h(\Sigma(2,3,11)) = -1$ is odd.
\end{example}

Let us try to adapt the construction of the Fr{\o}yshov invariant to the $\pin$-equivariant setting. In Section~\ref{sec:pin} we mentioned that $\swfh^{\pin}_*(Y; \F)$ is the homology of a complex composed of a triple of infinite $v$-tails (connected by the action of $q$) coming from the reducible, and a finitely generated piece coming from the irreducibles. In homology, we end up with three $v$-tails, which can end in various degrees. We obtain three different invariants $a(Y), b(Y), c(Y)$, given by the minimal degrees of nonzero elements in each of the tails. Since the tails have period $4$, and because the bottom element in the bottom $v$-tail in $\swfc^{\pin}_*(Y; \F)$ is in degree $-2n(Y, g)$, the invariants satisfy
$$ a(Y) \equiv b(Y)-1 \equiv c(Y)-2 \equiv -2n(Y,g) \equiv 2\mu(Y) \pmod 4.$$

Define
$$ \alpha(Y) = a(Y)/2, \ \ \beta(Y) = (b(Y)-1)/2, \ \ \gamma(Y) = (c(Y)-1)/2.$$
These are all $\Z$-lifts of the Rokhlin invariant. Moreover, they descend to $\Theta^H_3$, by a similar argument to the one for $h$. When we reverse the orientation of $Y$, by studying the long exact sequence \eqref{eq:tate} we find that the Tate homology (which has three $v$-tails that are infinite in both directions) is composed of the three $v$-tails of $\swfh^{\pin}_{*}(Y; \F)$ matched up with the three $v$-tails of $\swfh_{\pin}^{-*}(-Y; \F)$. The change in degree sign means that the bottom $v$-tail for $Y$ gets matched with the top $v$-tail for $-Y$, and vice versa, and the middle $v$-tails for $Y$ and $-Y$ get matched with each other. We deduce that:
$$ \alpha(-Y) = - \gamma(Y), \ \ \  \beta(-Y)=-\beta(Y), \  \ \gamma(-Y) = - \alpha(Y).$$ 
In conclusion, the middle invariant $\beta$ satisfies the desired properties (a) and (b), and this completes the proof of Theorem~\ref{thm:beta}.

\begin{example}
Consider $\Sigma(2,3,11)$ again. The three $v$-tails in the $\pin$-equivariant Floer homology \eqref{eq:pin11h} end in degrees $(a, b, c)=(4,1,2)$. We get that $(\alpha, \beta, \gamma)=(2,0,0)$, all even, in agreement with the fact that $\Sigma(2,3,11)$ has Rokhlin invariant zero.
\end{example}

To review, the key reason why $\beta$ worked better than the Fr{\o}yshov invariant for our purposes was that the cohomology of $B\pin$ is $4$-periodic, whereas the cohomology of $BS^1$ is $2$-periodic.

\subsection{Failure of additivity}
\label{sec:NoAdd}
Fr{\o}yshov proved in \cite{FroyshovSW} that his $h$-invariant gives a homomorphism from $\Theta^H_3$ to $\Z$. Let us briefly explain the curious facts that $h$ is a homomorphism and $\beta$ is not.

Fr{\o}yshov's proof starts with the observation that his construction works also for disjoint unions of homology spheres. One can define the group $\Theta^H_3$ as generated by all (possibly disconnected) manifolds with $H_1(Y; \Z)=0$, modulo homology cobordism. The connected sum of two manifolds is homology cobordant to their disjoint union, so we can use the latter to define addition in $\Theta^H_3$. Recall from Section~\ref{sec:du} that the Floer stable homotopy invariant of $Y_0 \amalg Y_1$ is obtained as the smash product of the invariants for the two factors. This implies that the $S^1$-equivariant Floer cochain complex for $Y_0 \amalg Y_1$ is obtained by tensoring those for $Y_0$ and $Y_1$ over the ground ring $H^*(BS^1; \f) = \f[U]$, where $\f$ is our chosen field. Since $\f[U]$
is a principal ideal domain, every $\f[U]$-chain complex is quasi-isomorphic to its homology. This implies that, in the derived category of $\f[U]$-modules, each $\swfc^{S^1}_*(Y_i)$ can be decomposed as a direct sum of a $U$-tail and a finite piece, and the two $U$-tails coming from $Y_0$ and $Y_1$ get tensored together to yield the tail for $Y_0 \amalg Y_1$. From here it is easy to deduce that $h(Y_0 \amalg Y_1) = h(Y_0) + h(Y_1)$.

The argument above fails in the $\pin$ case, because $H^*(B\pin; \F) = \F[q, v]/(q^3)$ is not a PID. In fact, one can give an explicit counterexample as follows. Let $Y = \Sigma(2,3,11)$, with Floer stable homotopy type $X = \swf(Y)$, the unreduced suspension of $\pin$. (Compare Example~\ref{ex:ex}.) Then $\swf(Y \amalg Y) = X \wedge X$. The smash product of two unreduced suspensions is the unreduced suspension of the join product. Moreover, the join of $\pin$ with itself is a $\pin$-bundle over the wedge sum $S^2 \vee S^2 \vee S^1$. Starting from these observations we can compute $\swfh_*^{\pin}(Y \amalg Y; \F)$ to be
 \[
 \xymatrixcolsep{.7pc}
\xymatrixrowsep{0pc}
\xymatrix{
  \F  & 0 & \F  & \F \ar@/_1pc/[l]_{q}  & \F \ar@/_1pc/[l]_{q} \ar@/^1pc/[llll]^{v} & 0 & \dots  \ar@/^1pc/[llll]^{v} & \dots \ar@/^1pc/[llll]^{v} & \dots \ar@/^1pc/[llll]^{v}\\
   & \oplus  & & & & & & &\\
    & \F^2 & & & & & & & 
} 
\]
with the leftmost element in degree $2$. This shows that the values of $(\alpha, \beta, \gamma)$ for $Y \amalg Y$ are $(2,2,0)$. Therefore, 
$$\beta(Y \amalg Y) = 2 \neq \beta(Y) + \beta(Y) = 0.$$

By considering disjoint unions of several copies of $\Sigma(2,3,11)$ (with both possible orientations), it can be shown that no non-trivial linear combination of $\alpha, \beta$ and $\gamma$ is additive.

\section{Other applications}
\label{sec:other}

\subsection{Smooth embeddings of four-manifolds with boundary} An important use of the Seiberg-Witten invariants in dimension four is the detection of exotic smooth structures. For example, one can prove  that the $K3$ surface has infinitely many distinct smooth structures. To show the similar result for the connected sum of several copies of the $K3$ surface, the usual Seiberg-Witten invariants do not suffice. However, Bauer and Furuta \cite{BauerFuruta} used finite dimensional approximation to define a stable homotopy version of the Seiberg-Witten invariant. This can be employed to show that the connected sum of up to four copies of $K3$ admits exotic smooth structures \cite{Bauer}.

The cobordism maps on Floer homology described in Section~\ref{sec:cobordism} produce relative Seiberg-Witten invariants of four-manifolds with boundary. A typical application of these is obstructing the embedding of a given four-manifold with boundary into a closed four-manifold. For example, recall that the $K3$ surface contains a nucleus $N(2)$ (the neighborhood of a cusp fiber and a section in an elliptic fibration), and $N(2)$ has boundary $-\Sigma(2,3,11)$. For $p, q > 0$ relatively prime with $(p, q) \neq (1,1)$, one can do logarithmic transformations of multiplicities $p$ and $q$ along elliptic fibers in $N(2)$ to obtain an exotic nucleus $N(2)_{p,q}$. Stipsicz and Szab\'o \cite{StSz} used Seiberg-Witten theory to show that the $K3$ surface cannot contain an embedded copy of the exotic nucleus of this form.

By using the stable homotopy version of the cobordism maps (also described in Section~\ref{sec:cobordism}), we obtain a similar obstruction for some connected sums:

\begin{theorem}[\cite{GluingBF}]
For any $p, q > 0$ relatively prime, with $(p, q) \neq (1,1)$, the exotic nucleus $N(2)_{p,q}$ cannot be smoothly embedded into the connected sum $K3 \# K3 \# K3$.
\end{theorem}

\subsection{Definite intersection forms of four-manifolds with boundary}
In Section~\ref{sec:beta} we explained that the Fr{\o}yshov invariant $h$ is an invariant of homology cobordism. More generally, if we have a negative-definite cobordism $W$ between integral homology spheres $Y_0$ and $Y_1$, similar methods produce the inequality
\begin{equation}
\label{eq:froyshov}
 h(Y_1) \geq h(Y_0) + \bigl(c_1(\t)^2  +  |\sigma(W)| \bigr)/8,
\end{equation}
where $\t$ can be any $\spinc$ structure on $W$; see \cite[Theorem 4]{FroyshovSW} or \cite[Theorem 39.1.4]{KMbook}. (A variant of this result first appeared in \cite{Froyshov1}.) The inequality \eqref{eq:froyshov} gives constraints on the possible intersection forms of negative-definite manifolds with fixed boundary. For example:

\begin{theorem}[Fr{\o}yshov \cite{Froyshov1}]
Let $W$ be a smooth, compact, oriented four-manifold with boundary the Poincar\'e sphere $P$. If the intersection form of $W$ is of the form $m\langle -1 \rangle  \oplus J$ with $J$ even and negative definite, then $J = 0$ or $J = -E_8$.
\end{theorem} 

The same results can be obtained by the Conley index method; see \cite{Spectrum, beta}. Further, the homology cobordism invariants $\alpha$, $\beta$ and $\gamma$ satisfy analogues of \eqref{eq:froyshov}, but they only apply to spin four-manifolds $(W, \t)$, with $c_1(\t)=0$.

\subsection{Indefinite intersection forms of four-manifolds with boundary}
For closed four-manifolds, Furuta \cite{Furuta} found the following constraint on the possible indefinite intersection forms: If $W$ is spin, then
\begin{equation}
\label{eq:108}
 b_2(W) \geq \frac{10}{8} |\sigma(W)| + 2.
\end{equation}
(The $11/8$ conjecture \cite{Matsumoto} claims the stronger inequality $b_2(W) \geq \tfrac{11}{8} |\sigma(W)|$.) Furuta's proof involved finite dimensional approximation for the Seiberg-Witten equations on $W$, which yields a $\pin$-equivariant map between representation spheres. The constraint \eqref{eq:108} is then obtained by studying the effect of this map on $\pin$-equivariant K-theory.

The theory of $\pin$-equivariant Seiberg-Witten Floer stable homotopy types, as presented in Section~\ref{sec:pin}, allows the generalization of \eqref{eq:108} to spin four-manifolds with boundary \cite{kg}. (Similar results were announced by Mikio Furuta and Tian-Jun Li.) Precisely, if $W$ has boundary an integral homology sphere $Y$, we get inequalities similar to \eqref{eq:108}, but involving a term that depends on the $\pin$-equivariant K-theory of $\swf(Y)$. 

For example, if $W$ has boundary a Brieskorn sphere of the form $\Sigma(2,3,12n+1)$, we get exactly the same inequality \eqref{eq:108}. In the cases $n=1$ or $2$, this result can be obtained more easily by observing that $\Sigma(2,3,13)$ and $\Sigma(2,3,25)$ are homology cobordant to $S^3$; therefore one can cap off $W$ with a homology ball, and apply Furuta's theorem to the resulting closed four-manifold. However, this simpler method does not work for larger $n$, since it is not known whether the Brieskorn spheres $\Sigma(2,3,12n+1)$ bound homology balls for $n \geq 3$.

\subsection{Covering spaces} Covers play a fundamental role in the topology of $3$-manifolds, but understanding their relationship with Floer homology is challenging.

Suppose we have a regular cover $\pi: \widetilde{Y} \to Y$ relating rational homology spheres, and let $\s$ be a $\spinc$ structure on $Y$.  The group $G$ of deck transformations of $\pi$ acts on the configuration space $\Conf(\widetilde{Y}, \pi^* \s)$, with fixed point set $\Conf(Y, \s)$. Introducing this additional $G$-symmetry into Floer theory is difficult by a Morse-theoretic approach. On the other hand, if we do finite dimensional approximation, we find that the Conley index $\tilde{I}^{\nu}$ corresponding to $\widetilde{Y}$ has a $G$-action, whose fixed point set is the Conley index $I^{\nu}$ for $Y$. By applying the classical Smith inequality \cite{Smith} to these Conley indices, we obtain the following:

\begin{theorem}[\cite{LidmanM}]\label{thm:Smith}
Suppose that $\widetilde{Y}$ and $Y$ are rational homology spheres, and $\pi:\widetilde{Y} \to Y$ is a $p^n$-sheeted regular covering, for $p$ prime. Let $\s$ be a $\spinc$ structure on $Y$. Then, the following inequality holds: 
\begin{equation}\label{eq:swfsmith}
\sum_i \dim \swfh_i(Y,\s); \f_p) \leq \sum_i \dim \swfh_i(\widetilde{Y},\pi^*\s; \f_p).
\end{equation}     
Here, $\f_p$ is the field with $p$ elements, and $\swfh_*$ denotes non-equivariant Seiberg-Witten Floer homology, i.e., the non-equivariant (reduced) singular homology of the Conley index, with the usual shift in degree by $\dim V^0_{-\nu} +2n(Y, \s, g)$.
\end{theorem}

\bibliographystyle{amsalpha}
\bibliography{biblio}

\end{document}